\documentclass[review]{elsarticle}

\usepackage{lineno,hyperref}
\modulolinenumbers[5]
\usepackage[nottoc]{tocbibind}

\usepackage{amsmath}
\usepackage{bm}

\begin{document}

\begin{frontmatter}

\title{A path integral approach to Bayesian inference \\ 
in  Markov  processes}

\author[mymainaddress1]{Toshiyuki Fujii}

\author[mymainaddress2]{Noriyuki Hatakenaka\corref{mycorrespondingauthor}}
\cortext[mycorrespondingauthor]{Corresponding author}
\ead{noriyuki@hiroshima-u.ac.jp}

\address[mymainaddress1]{Department of Physics, Asahikawa Medical University, Midorigaoka-higashi, Asahikawa, 078-8510, Japan}

\address[mymainaddress2]{Graduate School of Integrated Arts and Sciences, Hiroshima University, Higashi-Hiroshima 739-8521, Japan}

\begin{abstract}
We formulate Bayesian updates in Markov processes 
by means of path integral techniques and derive  the imaginary-time Schr\"{o}dinger equation with likelihood to direct the inference incorporated as a potential for the posterior probability distribution. 
\end{abstract}

\begin{keyword}
Bayesian inference  \sep Path integral technique \sep  Imaginary-time Schr\"{o}dinger equation
\MSC[2010] 
         62F15 %(Bayesian inference)
\sep 58D30 %(Feynman path integrals)
\sep 60H15 %(Stochastic partial differential equations) 
\end{keyword}

\end{frontmatter}

\linenumbers

%===============
\section{Introduction}
%===============
Statistical inference is a fundamental technology indispensable in modern society. In particular, Bayesian inference is explosive growth in a wide range of all kind of human activities, 
including philosophy, medicine, sport, law as well as science and technology \cite{o2010oxford, kruschke2017bayesian}. 
This uses a procedure to  update the  probability sequentially  based on the Bayes' theorem that incorporates new evidence into prior probabilities to form updated probability estimates. The form of this theorem is 
\begin{align}
P(H|E)=\frac{P(E|H)}{P(E)}P(H),
\end{align}
where $P(E|H)$ is the probability for the evidence $E$ given a hypothesis $H$, 
$P(H)$ is the prior probability that the hypothesis is true, 
$P(H|E)$ is the posterior probability that the hypothesis is true given the evidence, 
and $P(E)=\sum_HP(E|H)P(H)$ is the probability for the evidence averaged over all hypotheses.  In this way, 
the Bayes' theorem relates the prior probability of the hypothesis before getting the evidence $P(H)$ to the posterior probability of the hypothesis after getting the evidence $P(H|E)$ through 
the likelihood ratio $P(E|H)/P(E)$.

As is well known, the time evolution of the conditional probability (the transition probability) of the Markov process follows the Chapman-Kolmogorov equation  \cite{Chapman34,Kolmogoroff1931}; 
\begin{align}
p(x_{i+1},t_{i+1}| x_{i-1},t_{t-1})=\int p(x_{i+1},t_{i+1}| x_{i},t_{i})p(x_{i},t_{t}| x_{i-1},t_{i-1})dx_i,
\end{align}
where $p(x_{i}, t_{i}| x_j, t_j)$ is the conditional probability density function defined as the 
$p(x_{i}, t_{i}| x_j, t_j) dx_i dx_j \equiv P(x_{i} < X_{i} \leq x_{i} + dx_{i} |x_j < X_j \leq x_j + dx_j  )$
with $x_j$ and $X_j$ being the realization and the random variable at a certain time $t_j$. 
This expresses the chain condition of the Markov process that the state in the near future depends only on the state at the present time, not on past memories. This is consistent with the idea of infinitesimal time evolution in path integral techniques of quantum mechanics. In fact, Onsager and Machlup later incorporated Markov chain conditions into path integral techniques \cite{PhysRev.91.1505}. Along this line, this paper extends the Onsager-Machlup theory by incorporating the likelihood required by Bayesian inference into the Markov chain conditions.

%=========================================
\section{A path integral representation of Bayesian updates}
%=========================================

Suppose the stochastic processes under the probability distribution density functions $p_i(q_i)$ and $p_i(x_i)$ that are defined as $p_i(q_i)dq_i \equiv P(q_i < Q_i \leq q_i+ dq_i ) $ and  $p_i(x_i)dx_i \equiv P(x_i < X_i \leq x_i+ dx_i )$, respectively. Here the continuous random variables $Q_i$ and $X_i$ act as the roles of $H$ and $E$ in the Bayes' theorem
in the $i$-th update. 
Starting with the law of total probability 
\begin{equation}
p_{i+1 }(q_{i+1}) = \int d q_i p(q_{i+1}| q_i) p_i (q_i), 
\label{hatten}
\end{equation}
the first Bayesian update is expressed as 
\begin{equation}
p_1(q_1|x_1) =  \int d{q}_0 \frac{ p({x}_1| {q}_1 ) }{p({x}_1 )}  p(q_{1}| q_0) p_0(q_0).
\end{equation}
After repeating Bayesian updates of $N$ times, the posterior probability distribution density function is given as 
\begin{equation}
p_N( q_N | x_1, x_2, \cdots  x_N) = \int d q_0 K(q_N, q_0 ; x_1, x_2, \cdots x_N) p_0 (q_0),
\end{equation}
where the kernel $ K(q_N, q_0; x_1, x_2, \cdots x_N)$ is given as
\begin{align}
\lefteqn{
 K(q_N, q_0 ; x_1, x_2, \cdots x_N) } \notag \\
    &=\int dq_1 \int dq_2 \cdots \int dq_{N-1} \prod_{i=1}^{N} \frac{p(x_{i}| q_{i})}{p(x_{i})} p(q_i|q_{i-1}). 
 \label{karnel}
\end{align}
This describes repeated Bayesian updates representing the time evolution of stochastic processes starting from the prior probability distribution density function $p_0 (q_0)$ inferred from the repeated conditional probabilities.
Hereafter, we simply represent $p_i({q}_i | {x}_1, {x}_2, \cdots  {x}_{i-1}) = p_i ( q_i| \bm{x} )$ where $\bm{x}$ denotes the set of the realizations. 

Here, we introduce the parameter $t_i = i \epsilon$ with a small quantity $\epsilon$ and an integer $i$, in order to prepare for converting  discrete time to continuous time such as $q_i = q (t_i)$.  
In this notation, for example, the probability distribution density function and the transition probability distribution density function are expressed as 
$p_i( q_i =q | \bm{x} ) = p(q, t_i |  \bm{x} )$ and $p( q_j = q' | q_i = q) = p(q', t_j |q,  t_i)$, respectively. 
Let us consider the infinitesimal time evolution of $p(q', t_j |q,  t_i)$ from time $t_i$ to $t_j=t_{i+1}$. 
Since $
p(q', t_{i+1} |q,  t_i)= p(q', t_i+\epsilon |q,  t_i) $, then
\begin{eqnarray}
p(q', t_i+\epsilon |q,  t_i) 
&=& 
p(q', t_i |q,  t_i) 
+\left.\epsilon\frac{\partial p(q', t'_i |q,  t_i)}{\partial t'_i} \right|_{t'_i=t_i}
+ O(\epsilon^2) .
\label{p-expand}
\end{eqnarray}
Suppose that the stochastic process that we are assuming is a Wiener process. 
If the stochastic process is a Markov process, the random variable $q$ obeys the Ito-type stochastic partial differential equation and the probability distribution density function based on this random variable follows the Fokker-Planck (FP) equation. Furthermore, the FP equation results in a diffusion equation in the case of a Wiener process leading to Brownian motion which has no driving force. Therefore, the second term on the right hand side of Eq. \eqref{p-expand} is replaced by the diffusion equation; 
\begin{align}
\frac{\partial p(q', t' |q,  t)}{\partial t'}
=
\frac{1}{2}D\frac{\partial^2 p(q', t' |q,  t)}{\partial q'^2},
\end{align}
where $D$ is the diffusion constant. The resulting infinitesimal time-evolution equation is 
\begin{align}
p(q', t_i+\epsilon |q,  t_i) 
=
\left.\left( 1+\frac{\epsilon}{2} D \frac{\partial^2 }{\partial q'^2} \right)
 p(q', t'_i |q,  t_i) \right|_{t'_i=t_i}
+ O(\epsilon^2).
\end{align}
By using the relation $p(q', t'_i |q,  t_i)|_{t'_i=t_i}=\delta(q'-q)$, this is reduced to 
\begin{align}
p(q', t_i+\epsilon |q,  t_i) 
&=
\left( 1+\frac{\epsilon}{2}D\frac{\partial^2 }{\partial q'^2} \right)
\delta(q'-q)
+ O(\epsilon^2)  \notag \\
&=
\left( 1+\frac{\epsilon}{2}D\frac{\partial^2 }{\partial q'^2} \right)
\left(\int \frac{ds}{2\pi} e^{i(q'-q)s}\right)
+ O(\epsilon^2)  \notag \\
&=
\int \frac{ds}{2\pi}  \left( 1+\frac{\epsilon}{2}D(-s^2) \right)
e^{i(q'-q)s}
+ O(\epsilon^2)  \notag \\
&=
\int \frac{ds}{2\pi}  e^{-\frac{\epsilon}{2}Ds^2}
e^{i(q'-q)s}
+ O(\epsilon^2)  \notag \\
&=
\frac{1}{\sqrt{2\pi D\epsilon}}e^{-\frac{(q'-q)^2}{2D\epsilon}}
+ O(\epsilon^2).  
\end{align}
Here we define the velocity of the random variable by $\dot{q}=(q'-q)/\epsilon$.
The equation is then given as 
\begin{align}
p(q', t_i+\epsilon |q,  t_i) 
&=
\frac{1}{\sqrt{2\pi D\epsilon}}e^{-\frac{\dot{q}^2}{2D}\epsilon}
+ O(\epsilon^2).  
\end{align}
This is consistent with the Onsager-Machlup formula  \cite{PhysRev.91.1505} in the path integral representation of stochastic processes. 

Now, let us incorporate the likelihood which is the most significant factor in Bayesian update into the path integral form. 
We adopt the exponential function used in the previous studies  \cite{Bialek96fieldtheories} as likelihood $ p(x, t_i| q, t_i )$; 
\begin{equation}
p(x, t_i| q, t_i  ) = M \exp
[
-V (q(t_i)-x(t_i)) \epsilon 
],
\end{equation}
where $M$ is a normalization constant. 
This form satisfies the requirement for the probability density function that it has a non-negative value in the range from $-\infty$ to $\infty$. Thus, 
the kernel is represented as 
\begin{align}
\lefteqn{
 K(q(t_N), q(t_0), \{x(t_i), i= 1, 2, \cdots N \} ) 
 }\nonumber \\ 
    &= M^N \int dq(t_1) \int dq(t_2) \cdots \int dq(t_{N-1}) \notag \\
& \hspace{1cm}\times   \exp\left[{-\epsilon  \sum_{i=1}^{N} \left\{ \frac{\dot q(t_n)^2}{2D}+V (q(t_n) -x(t_n)) \right\}  }\right].
    \label{risank}
\end{align}
Now let us take a limit $\epsilon \rightarrow 0$ under the condition
$t_N = N \epsilon \equiv t' $. 
The discrete variable $ t_i = i \epsilon $ becomes dense  at time interval $ [0, t'] $. 
The set of realizations in stochastic processes  $ \bm{q}$ and $\bm{x}$ can be represented by continuous functions $ q (t) $ and $ x (t) $.
Thus, the probability distribution density function $p_i(q_i )$ is also expressed as $p(q, t )$. 
In this limit. Eq. \eqref{risank} is just a path integral description of Bayesian updates; 
\begin{equation}
 K(q(t'), q(t ), {x(t)}) =\int \mathcal{D}[q]   e^{-  S[q(t), \dot q(t), x(t)] },
 \label{contkarnel}
\end{equation}
with
\begin{equation}
 S[q(t), \dot q(t), x(t)]= \int^{t'}_{0} \left\{\frac{\dot q(t)^2}{2D} +V (q(t)-x(t)) \right\}  dt, 
 \label{S}
\end{equation}
where 
$\int \mathcal{D}[q]  \equiv M^N \int dq(t_1) \int dq(t_2) \cdots \int dq(t_{N-1})
$. 
This action $S$ describes a particle moving in the random potential $V$. 
This is a main result of this paper. The likelihood, which is the most important factor in Bayesian update, is incorporated into the path integral form as a potential. This result seems to be quite natural. As is well known, the potential leads the force, which determines the direction of particle movement in mechanics. It is exactly consistent with the role of likelihood in Bayesian update.

\section{An imaginary-time Sch\"{o}dinger equation}
Let us derive the dynamical equation of the posterior probability distribution density associated with the Bayesian update based on the path integral formula for Bayesian updates derived in the previous section. 
Under an infinitesimal interval $\epsilon$ from $t$ to $t'$ ( 
$t' = t+ \epsilon$), 
\begin{eqnarray}
p(q', t' |\bm{x} ) = M \int dq  e^{-\epsilon \left[  \frac{1}{2D}\left\{ \frac{q' - q}{ \epsilon} \right\}^2 + V(q-x(t)) \right]} p(q, t |\bm{x}).
\label{baysean2}
\end{eqnarray}
Let us take the Taylor series expansions for probability distribution density functions on both sides of the equation. The probability distribution density function $p(q', t' |\bm{x} )$ is expanded with respect to $\epsilon$ as
\begin{align*}
p(q', t+\epsilon|\bm{x} )  = p(q', t  |\bm{x} ) + \frac{\partial p(q', t|\bm{x} ) }{\partial t} \epsilon + O(\epsilon^2).  
\end{align*}
On the other hand, the probability distribution density function 
$ p(q, t |\bm{x} )$ on the right hand side can be expressed as $p(q' -\eta, t |\bm{x} )$ to meet the variable on the left side by introducing a variable $\eta =q'-q $, and is expanded with respect to $\eta$, 
\begin{align*}
\lefteqn{
p(q'-\eta, t |\bm{x} )}\notag \\
&  = p(q', t|\bm{x} ) - \frac{ \partial p(q',t|\bm{x} )}{\partial q'} \eta + \frac{1}{2}  \frac{ \partial^2 p(q',t|\bm{x} )}{\partial q'^2} \eta^2  + O(\eta^3).   
\end{align*}
Then, we obtain 
\begin{align}
\lefteqn{
p(q', t|\bm{x} ) + \frac{ \partial p(q', t|\bm{x} )}{\partial t} \epsilon} \notag \\
&= M \int d\eta  e^{- \frac{\eta^2}{2\epsilon}} \left[ 1 - \epsilon  V(q-x(t))  \right]\notag \\
&\hspace{1cm}\times
\left[  p(q', t |\bm{x}  ) + \frac{ \partial p(q', t|\bm{x} )}{\partial q'} \eta + \frac{1}{2}  \frac{ \partial^2 p(q', t|\bm{x} )}{\partial q'^2} \eta^2   \right].
\label{baysean3}
\end{align}
After Gaussian integration with 
\begin{align*}
\int_{-\infty}^{\infty} \eta^{2n} \exp (-\eta^2/2D\epsilon )  d\eta= (2n-1)!!(D \epsilon)^n\sqrt{2\pi D \epsilon}, 
\end{align*}
we get   
\begin{align}
\lefteqn{
p(q', t|\bm{x} ) + \frac{ \partial p(q', t|\bm{x} )}{\partial t} \epsilon
}\notag \\
&= M   \sqrt{2\pi D \epsilon} \left[   1   - \epsilon  V(q-x(t))  \right]\left[  p(q',  t |\bm{x}  )  + \frac{D}{2}  \frac{ \partial^2 P(q', t|\bm{x} )}{\partial q'^2} \epsilon  \right] \notag  \\
&=
M   \sqrt{2\pi D \epsilon} p(q', t |\bm{x}  ) \notag \\
& \hspace{1cm}+ 
M   \sqrt{2\pi D \epsilon^3}  
\left\{
 -   V(q-x(t)) 
+    \frac{D}{2}  \frac{ \partial^2 p(q', t|\bm{x} )}{\partial q'^2} 
 \right\}
+ O(\epsilon^2). 
\label{baysean4}
\end{align}
Comparing the first term on both sides of equation yields 
$M= (2\pi D \epsilon )^{-1/2}$. 
The remaining terms eventually form 
\begin{equation}
- \frac{\partial p(q, t|\bm{x} ) }{\partial t} = \left\{ - \frac{D}{2} \frac{\partial^2}{\partial q^2} + V(q-x(t)) \right\} p(q, t|\bm{x} ). 
\label{iseq}
\end{equation}
This describes the time  evolution of the posterior probability density function in Bayesian inference in Markov processes. 
This is nothing but the Sch\"odinger equation in quantum mechanics except that it is imaginary time. 
This equation is very advantageous for Bayesian inference in a way different from the path integral form because there are many technological accumulations in the  Schr\"odinger equation along with the development of quantum mechanics. In the next section, we will look at some of these advantages.

\subsection{Example: Gaussian likelihood}
Let us consider the time evolution of the probability distribution density function using the imaginary-time Schr\"odinger equation for Bayesian update obtained in the previous section by adopting the Gaussian likelihood as a concrete example given by 
\begin{equation}
V= \frac{[q-x(t)]^2}{2\Delta^2}. 
\end{equation}
The imaginary-time Schr\"odinger equation for Bayesian update to be solved is then  
\begin{equation}
- \frac{\partial p(q, t |\bm{x}   ) }{\partial t} = \left\{ - \frac{D}{2} \frac{\partial^2}{\partial q^2} +  \frac{[q-x(t)]^2}{2\Delta^2}  \right\} p(q, t |\bm{x}   ). 
\label{iseqho}
\end{equation}
This equation is the imaginary time version of the Schr\"odinger equations of forced harmonic oscillators with a random force $x(t)$
\begin{align}
F=-\frac{dV}{dq}=-\frac{q}{\Delta^2}q+\frac{q}{\Delta^2}x(t). 
\end{align}
This partial differential equation can simply be solved using the Feynman-Kac formula \cite{kac1949distributions, cheng1983exact}.  
However, we will seek the solution using the operator method developed in quantum mechanics to show the usefulness of the imaginary-time Schr\"odinger equation in Bayesian inference. 

By introducing  new operators called creation ($ a$) and annihilation ($ a^{\dagger}$) operators in quantum mechanics; 
\begin{align}
q &= \sqrt{\frac{\Delta \sqrt{D} }{2}} ( a+   a^{\dagger}), \\
\frac{\partial }{\partial q} 
&= \sqrt{\frac{1 }{2\Delta \sqrt{D}  }} ( a-  a^{\dagger}), 
\end{align}
with $[ a,  a^\dagger ] \equiv  a  a^\dagger -   a^\dagger  a  =1$,  
the equation is expressed by 
\begin{equation}
- \frac{\partial p(q, t|\bm{x} ) }{\partial t} = \left\{  \frac{\sqrt{D} }{\Delta } \left(  a^\dagger  a  + 1/2 \right) + \frac{x(t)^2}{2\Delta^2}  + x (t) \sqrt{\frac{\sqrt{D}  }{2 \Delta^3 }} ( a +   a^\dagger)    \right\} p(q, t |\bm{x}   ). 
\label{iseqhorhs}
\end{equation}
The first two terms on the right hand side represent operators corresponding to Hamiltonians for a simple harmonic oscillator and a random force in quantum mechanics.
On the other hand, the third term expresses their interaction. 
We focus on the time evolution of the probability distribution density function 
$p'(q, t |\bm{x} )$
due to the interaction using the operator 
\begin{equation}
 U_0 (t) = \exp{ \left\{  \frac{-\sqrt{D} }{\Delta } \left(  a^\dagger  a  + \frac{1}{2} \right) t  - \int^{t}_0 \frac{x(t' )^2}{2\Delta^2} dt'  \right\} },
\end{equation}
associated with $p(q, t |\bm{x}   ) = U_0 p'(q, t |\bm{x} )$, 
this reads 
\begin{equation}
- \frac{\partial p'(q, t |\bm{x} ) }{\partial t} =  x (t)  \sqrt{\frac{\sqrt{D}  }{2 \Delta^3 }} U_0(-t) ( a +   a^\dagger)  U_0(t) p'(q, t|\bm{x}).
\label{seqint}
\end{equation}
According to the Baker-Campbell-Hausdorff formula  \cite{baker1905alternants} for any two operators, say $A$ and $B$,  
\begin{equation}
e^A B e^{-A} = B + [A, B] + \frac{1}{2} [A, [A, B]] + \cdots \frac{1}{n!} [ A, [A, \cdots [A, B] \cdots] + \cdots ,
\end{equation}
the terms $U_0(-t) (\hat a +  \hat a^\dagger)  U_0(t)$ are expressed as 
\begin{equation}
U_0(-t)  (a +a^\dagger )  U_0(t) = e^{-\frac{\sqrt{D} t}{\Delta}} a + e^{\frac{\sqrt{D} t}{\Delta}} a^\dagger .
\end{equation}
Thus, Eq. \eqref{seqint} is now represented as
\begin{align}
- \frac{\partial p'(q, t |\bm{x}) }{\partial t}  & =  x (t)  \sqrt{\frac{\sqrt{D}  }{2 \Delta^3 }} \left( e^{-\frac{\sqrt{D} t}{\Delta}} a +  e^{\frac{\sqrt{D} t}{\Delta}} a^\dagger \right)  p'(q, t |\bm{x}   ) 
\notag \\
&\equiv H(t) p'(q, t |\bm{x}). 
\label{seqint2}
\end{align} 
The formal solution of this equation is well known and given as 
\begin{equation}
 p'(q, t |\bm{x}) = T\left[  \exp \left\{ -\int_0^{t} H(t') dt'  \right\} \right] p'(q, 0),
\end{equation}
where $T$ is the time ordering operator 
\begin{align}
T[A(t), B(t' )] 
= 
\begin{cases}
A(t)B (t') & \text{for } t >  t' \\
B (t')   A(t) & \text{for } t<  t' . 
\end{cases}
\end{align}
This is also expressed by \cite{9780471887027}
\begin{align}
 p'(q, t |\bm{x} )
 =
 \exp \left[ - \int^{t}_0 H(t') dt'  - \frac{1}{2}  \int^{t}_0 dt'  \int^{t'}_0dt'' \left[ H(t'),   H(t'') \right]   \right] p(q, 0  ),
 \label{tint}
\end{align}
since $[H(t ),  [H(t'), H(t'') ]] = 0$ from 
\begin{align}
\left[ H(t'),   H(t'') \right] &=  x (t')   x (t'') \frac{\sqrt{D} }{2 \Delta^3 }  \left[ ( e^{-\frac{\sqrt{D} t'}{\Delta}} a +  e^{\frac{\sqrt{D} t'}{\Delta}} a^\dagger ),    ( e^{-\frac{\sqrt{D}t''}{\Delta}} a +  e^{\frac{\sqrt{D}t''}{\Delta}} a^\dagger ) \right] \notag \\
&=
  x (t')   x (t'') \frac{\sqrt{D}  }{2 \Delta^3 }  \left[ ( e^{-\frac{\sqrt{D} (t'-t'')}{\Delta}} [ a,  a^\dagger] + e^{\frac{\sqrt{D} (t'-t'')}{\Delta}} [   a^\dagger, a] \right] \notag \\
  &=
 - x (t')   x (t'') \frac{\sqrt{D} }{ \Delta^3 }   \sinh \left( \frac{\sqrt{D}(t'-t'')}{\Delta} \right) . 
 \label{koukan}
\end{align}
The time evolution of the probability distribution density function which focused on the interaction is then given as 
\begin{align}
 p'(q, t |\bm{x}) 
&= e^{f(t)}
 \exp{\left[ -  \int_0^t dt'  x (t') \sqrt{\frac{\sqrt{D} }{2\Delta^3 }}  \left( ae^{- \frac{\sqrt{D} t'}{\Delta } }+  a^\dagger  e^{\frac{\sqrt{D} t'}{\Delta } } \right)\right]}  p'(q, 0), 
\end{align}
where $f(t) \equiv ( {\sqrt{D}}/ {2\Delta^3 } )\int_0^{t}dt' \int_0^{t'} dt''  x(t')  x(t'')  \sinh ( {\sqrt{D}(t'-t'')}/{\Delta} ) $. 
This can be reverted to the posterior probability distribution density function using the relation 
$p(q, t|\bm{x}) = U_0(t) p'(q, t |\bm{x}  )$, resulting in   
\begin{align}
 p(q, t |\bm{x}   ) 
 &= e^{f(t)}
 \exp \left[ -  \int_0^t dt'  x (t') \sqrt{\frac{1 }{2\Delta^3 }}  U_0(t) ( ae^{- \frac{\sqrt{D} t'}{\Delta } }+  a^\dagger  e^{\frac{ \sqrt{D}t'}{\Delta } } ) 
 U_0^{-1 } \right] \notag \\
 & \hspace{1cm} \times U_0(t)p(q, 0) \notag \\
  &= e^{f(t)}
 \exp \left[ -  \int_0^t dt' x (t') \sqrt{\frac{\sqrt{D} }{2\Delta^3 }}   ( ae^{- \frac{\sqrt{D}( t'-t)}{\Delta } }+  a^\dagger  e^{\frac{ \sqrt{D}(t'-t)}{\Delta } } ) \right] \notag \\
 &\hspace{1cm}\times U_0(t)p(q, 0), 
 \label{kai1}
\end{align}
where the relations  $p(q, 0) = p'(q, 0)$ and $U_0(t) F(A) U_0^{-1}(t) = F(U_0(t)A U_0^{-1})$
that holds for any function $F(A)$ of the operator $A$ are used. 

The prior probability distribution density function density is assumed to be given by the superposition of the harmonic oscillator solutions $f_n$ as follows
\begin{equation}
p(q, 0)	= \sum_n c_n f_n(q),  
\end{equation}
where 
\begin{equation}
f_n (q) = \frac{1}{\left({\pi \Delta \sqrt{D}}\right)^{1/4} \sqrt{2^n n!}  } e^{-\frac{q^2}{2\Delta \sqrt{D}}} H_n \left( \frac{q}{\sqrt{\Delta\sqrt{D} }}\right),
\end{equation}
with $H_n(x)$ being the $n$-th order Hermite polynomial at $x$. 
The expansion coefficient $c_n$ is given by 
$c_n = \int dq f_n(q) P_0 (q, 0)$. 
Starting with this prior probability distribution density function, the posterior probability distribution density function is given as follows; 
\begin{align}
 p(q, t|\bm{x}) 
  &= e^{f(t)}
 \exp \left[ -  \int_0^t dt'  x (t') \sqrt{\frac{\sqrt{D} }{2\Delta^3 }}   \left( ae^{- \frac{\sqrt{D}( t'-t)}{\Delta } }+  a^\dagger  e^{\frac{ \sqrt{D}(t'-t)}{\Delta } } \right) \right] 
 \notag \\
 &\hspace{1.0cm} \times  U_0(t)\left(\sum_{n=0}^{\infty} c_n f_n (q)\right) \notag \\
 &=
     e^{f(t) - \int^{t}_0 \frac{x(t' )^2}{2\Delta^2} dt' }
    \sum_{n=0}^{\infty} c_n  
  e^{  \frac{-\sqrt{D}}{\Delta } \left( n  + \frac{1}{2} \right) t  }  \notag \\
 & \hspace{0.5cm}
  \times \exp \left[ -  \int_0^t dt'  x (t') \sqrt{\frac{\sqrt{D} }{2\Delta^3 }}   \left( ae^{- \frac{\sqrt{D}( t'-t)}{\Delta } }+  a^\dagger  e^{\frac{ \sqrt{D}(t'-t)}{\Delta } } \right) \right]f_n (q).  
 \label{kai2}
\end{align}
Note that 
\begin{align}
e^{  -  \int_0^t dt'  x (t') \sqrt{\frac{\sqrt{D} }{2\Delta^3 }}   \left( ae^{- \frac{\sqrt{D}( t'-t)}{\Delta } }+  a^\dagger  e^{\frac{\sqrt{D}( t'-t)}{\Delta } } \right)  
}
=
e^{ -   x''(t)q}
e^{   -  x' (t)   \frac{\partial }{\partial q} }
e^{ g(t) },
\label{cohe}
\end{align}
where 
\begin{align}
 x'(t) &=  \frac{\sqrt{D}}{\Delta}   \int_0^t dt'  x (t')    \sinh \left\{ { \sqrt{D}(t-t')} / {\Delta }  \right\}, \\
 x'' (t)&= { \frac{1}{\Delta^2} }  \int_0^t dt'  x (t')    \cosh \left\{\frac{\sqrt{D}( t-t')}{\Delta }  \right\}, \\
 g(t) &=  { \frac{\sqrt{D}}{2\Delta^3} }  \int_0^t dt'  \int_0^t dt''  x (t')   x (t'')  
\notag \\
&\hspace{1cm}\times \cosh \left\{\frac{ \sqrt{D}(t-t')}{\Delta }  \right\}  \sinh \left\{\frac{\sqrt{D}( t-t'')}{\Delta }  \right\},    
 \end{align}
and  
$
e^{c \frac{\partial }{\partial q}} f(q) = f(q+c). $

As a result, we reach the final expression for the posterior probability distribution density function $p(q, t |\bm{x}  )$ as 
\begin{align}
 p(q, t |\bm{x}) 
    &=
     e^{f(t) +g(t) - \int^{t}_0 \frac{x(t' )^2}{2\Delta^2} dt' } 
   e^{ - x''(t) q   } \notag \\
 &\hspace{1cm}\times \sum_{n=0}^{\infty} c_n  e^{ -  \frac{\sqrt{D}}{\Delta } \left( n  + \frac{1}{2} \right) t  }
  f_n (q - x'(t) ). 
 \label{kai3}
\end{align}
We can find the role of Bayesian update in stochastic process in this solution obtained from stationary solution $f_n$ of the Schr\"odinger equation. 
Evidence incorporated in Bayesian update is embedded in the argument of the function $f_n$. 
This means that the center of the posterior distribution density function is shifted with every update. 
In addition, 
the first $n$-independent exponent is interpreted as the global phase corresponding to the geometric phase. 
These features are difficult to extract in the path integral form like Feynman-Kac form. 
Therefore, this is the greatest advantage of analysis that makes it possible to elucidate elementary processes of stochastic processes using the Schr\"odinger equation.

%================
\section{Summary}
%================
We have formulated the Bayesian inference in  Markov processes by means of path integral techniques 
and have succeeded in incorporating the likelihood which is the key in Bayesian update as a potential into the path integral form. 
Then we have derived the imaginary-time Schr\"{o}dinger  equation from path integral representation of the evolution of the Bayesian inference in  Markov processes. 
Through an example of Gaussian likelihood, we have showed that the operator technique for the Schr\"odinger equation is effective for elucidating the elementary processes of stochastic processes.

%\section*{References}

\bibliography{mybibfile}

\end{document}